\documentclass[11pt]{amsart}
\newcommand{\BZ}{\mathbb{Z}}
\newcommand{\BN}{\mathbb{N}}
\newcommand{\BR}{\mathbb{R}}
\newcommand{\BT}{\mathbb{T}}

\newcommand{\BC}{\mathbb{C}}

\newcommand{\BP}{\mathbb{P}}

\newcommand{\caB}{\mathcal{B}}
\newcommand{\caS}{\mathcal{S}}

\newcommand{\caD}{\mathcal{D}}

\newcommand{\caM}{\mathcal{M}}

\newcommand{\caG}{\mathcal{G}}

\newcommand{\gC}{\Gamma}

\newcommand{\gL}{\Lambda}

\newcommand{\ga}{\alpha}

\newcommand{\gd}{\delta}

\newcommand{\gve}{\varepsilon}
\newcommand{\gc}{\gamma}
\newcommand{\gf}{\varphi}
\newcommand{\gl}{\lambda}
\newcommand{\gs}{\sigma}

\newcommand{\gx}{\chi}
\newcommand{\gp}{\psi}
\newcommand{\gm}{\mu}

\newcommand{\gz}{\zeta}

\newcommand{\bt}{\bar{\triangle}}

\newcommand{\ti}[1]{\tilde{#1}}
\newcommand{\p}{\prod}

\newcommand{\te}{\text}

\newcommand{\hra}{\hookrightarrow}
\newcommand{\ra}{\rightarrow}

\newcommand{\lra}{\longrightarrow}

\newcommand{\df}{\stackrel{\te{def}}{=}}

\newcommand{\ip}[2]{\left<{#1},{#2}\right>}
\newcommand{\Lim}[1]{\lim_{{#1} \ra \infty}}
\newcommand{\veC}[2]{({#1}_1,\ldots,{#1}_{{#2}})}
\newcommand{\Sum}[2]{\sum_{{#1}=1}^{{#2}}}
\newcommand{\Avr}[2]{\frac{1}{{#1}}\sum_{{#2}=1}^{{#1}}}

\def\lfa#1#2{\mathord{\hbox{\vtop{\kern-5pt\hbox{$#2$}}$\setminus$$#1$}}}

\swapnumbers
\theoremstyle{plain}
\newtheorem{lma}{Lemma}[section]
\newtheorem{thm}[lma]{Theorem}
\newtheorem{pro}[lma]{Proposition}

\newtheorem{cor}[lma]{Corollary}

\theoremstyle{definition}
\newtheorem{dfn}[lma]{Definition}
\newtheorem{rmr}[lma]{Remark}

\newtheorem{dsc}[lma]{}

\begin{document}

\title[Configurations in the Plane]{Nilfactors of $\BR^m$-actions and
configurations in sets of positive upper density in $\BR^m$} 
\author{T. Ziegler}
\maketitle

\begin{abstract}
We use ergodic theoretic tools to solve a
classical problem in geometric Ramsey theory.
Let $E$ be a measurable subset of $\BR^m$, with $\bar{D}(E)>0$.
Let $V=\{0,v_1,\ldots,v_{k}\} \subset \BR^m$. We show that for $r$ large
enough, we can find an isometric copy of $rV$ arbitrarily close
to $E$. This is a generalization of a theorem of Furstenberg, 
Katznelson and Weiss ~\cite{FuKaW} showing a similar property for 
$m=k=2$.  
\end{abstract}

\section{Introduction}

Let $E$ be a measurable subset of $\BR^m$.
We set
\[ 
\bar{D}(E):= \limsup_{l(S) \rightarrow \infty} \frac{m(S \cap E)}{m(S)},
\]
where $S$ ranges over all cubes in $\BR^m$, and $l(S)$ denotes 
the length of a side of $S$.
$\bar{D}(E)$ is the upper density of $E$.
We are interested in configurations which are necessarily contained in $E$.
Furstenberg, Katznelson, and Weiss ~\cite{FuKaW} showed, using methods from ergodic theory, that if $E \subset \BR^2$, with $\bar{D}(E)>0$, all large distances 
in $E$ are attained. More precisely:
\begin{thm}[FuKaW]\label{T:Distances}
If $E \subset \BR^2$ with $\bar{D}(E)>0$, there exists $l_0$ such that for any
$l>l_0$ one can find a pair of points $x,y \in E$ with $\|x-y\|=l$.
\end{thm}
This result was also proved, using different methods, by Bourgain ~\cite{Bo}, and by Falconer and 
Marstrand ~\cite{FaMa}.
It is natural to ask if the same is valid for larger configurations.
Bourgain has shown by an example that this can not be done ~\cite{Bo}.

As some configurations may not be found in the set itself,
we try to find the configurations
arbitrarily close to the set.
In the same paper  Furstenberg, Katznelson, and Weiss ~\cite{FuKaW} 
show that with this weaker condition, 
one can find triangles in the plane:
\begin{thm}[FuKaW]
Let $E \subset \BR^2$ with $\bar{D}(E)>0$, and let $E_{\gd}$ denote
the points at distance $<\gd$ from $E$. Let $v,u \in \BR^2$, then 
there exists $l_0$ such that for $l>l_0$ and any $\gd>0$ there exists a triple 
$(x,y,z) \subset E_{\gd}^3$ forming a triangle congruent to $(0,lu,lv)$.
\end{thm}
The idea of the proof is to translate the geometric problem to a dynamical
problem, where $E$ corresponds to some measurable set $\ti{E}$, with
positive measure, in a measure preserving system $(X^0,\caB,\mu,\BR^2)$. 
The statement that $E_{\gd}$ contains a certain configuration, 
corresponds to a 
recurrence condition on the set $\ti{E}$. In the case of triangles
(configurations formed by $2$ vectors), the recurrence phenomenon in question
is reduced to the 
case where $(X^0,\caB,\mu,\BR^2)$ is a Kronecker action.
The problem for a general configuration reduces to the study of 
pro-nilsystems (defined later). We prove the following theorem: 
%
% MAIN THEOREM
%
\begin{thm}\label{thm:main}  
Let $E \subset \BR^m$ have positive upper density, 
and let $E_{\gd}$ denote the points of distance $< \gd$ from $E$. 
Let $(u_1,\ldots,u_k) \subset (\BR^m)^k$. 
Then there exists $l_0$ such that for any $l > l_0$, and any $\gd>0$ 
there exists $\{x_1,x_2,\ldots,x_{k+1}\} \in E_{\gd}^{k+1}$
forming a configuration congruent to $\{0,lu_1,\ldots,lu_{k}\}$.
\end{thm}

{\bf Acknowledgment}
I thank Sasha Leibman and Hillel Furstenberg for helpful comments.

%%%%%%%%%%%%%%%%%%%%%%%%%%%%%%%%%%%%%%%%%%%%%%%%%%%%%%
%%%%%%%%%%%%%%%%%%%%%%%%%%%%%%%%%%%%%%%%%%%%%%%%%%%%%%
\section{Translation of the Geometric Problem \\    %% 
         to a Dynamical Problem.}                   %%
%%%%%%%%%%%%%%%%%%%%%%%%%%%%%%%%%%%%%%%%%%%%%%%%%%%%%%
%%%%%%%%%%%%%%%%%%%%%%%%%%%%%%%%%%%%%%%%%%%%%%%%%%%%%%

We start by translating the geometric problem
to a dynamical problem. The translation as shown here was done
in ~\cite{FuKaW}. We bring it here for the sake of completeness.\vspace{.25in}

Let $E \subset \BR^m$, such that $\bar{D}(E)>0$. Define
\[
\gf(u):=min \{1,dist(u,E)\}.
\]
The functions $\gf_v(u)=\gf(u+v)$ form an equicontinuous, uniformly 
bounded family, and thus have compact closure in the topology of 
uniform convergence over bounded sets in $\BR^m$. Denote this 
closure by $X^0$. $\BR^m$ acts on $X^0$ by $T_v\gp(u)=\gp(u+v)$ for
$\gp \in X^0, u,v\in \BR^m$. $X^0$ is a compact metrizable space and
we can identify Borel measures on $X^0$ with functionals on $C(X^0)$.  
Since $\bar{D}(E)>0$, there exists a sequence of cubes $S_n$
such that
\[ 
\frac{m(S_n \cap E)}{m(S_n)} \lra \bar{D}(E)>0.
\]
We define a probability measure $\gm$ on $X^0$ as follows. 
We define the following probability measures:  
for $f \in C(X^0)$, let
\[
  \gm_n(f)=\frac{1}{m(S_n)}\int_{S_n} f(T_v\gf) dm(v)
\] 
We have for some subsequence $\{n_k\}$ 
\[
  \gm_{n_k} \overset{w*}{\lra} \gm.
\]
Set $f_0(\gp)=\gp(0)$, then $f_0$ is a continuous function on $X^0$.
We define $\tilde{E} \subset X^0$ by
\[
  \gp \in \tilde{E} \iff f_0(\gp)=0 \iff \gp(0)=0.
\] 
$\tilde{E}$ is a closed subset of $X^0$ and we have:
\[
  \mu(\tilde{E})= \lim_{l \ra \infty} \int_X (1-f_0(\gp))^l d\mu(\gp).
\]
%
% The measure of the set E is > 0
%
\begin{lma}
$\mu(\tilde{E})>0$.
\end{lma}

\begin{proof}
It suffices to show that for any $l$,
\[
 \int_X (1-f_0(\gp))^l d\mu(\gp) \ge \bar{D}(E).
\] 
Indeed
\begin{equation*}
  \begin{split}
    \int_X (1-f_0(\gp))^l d\mu(\gp)  
        &= \lim_{k \ra \infty}  \frac{1}{m(S_{n_k})}\int_{S_{n_k}} 
           (1-f_0(T_v\gf))^l dm(v)                                \\
        &= \lim_{k \ra \infty}  \frac{1}{m(S_{n_k})}\int_{S_{n_k}}
           (1-\gf(v))^l dm(v)                                     \\
        &\ge  \lim_{k \ra \infty} \frac{m(S_{n_k} \cap E)}{m(S_{n_k})} 
         = \bar{D}(E) > 0,
  \end{split}
\end{equation*}
since $\gf(v) = 0$ for $v \in E$.
\end{proof} 

The next proposition establishes the correspondence between $E$ and 
$\tilde{E}$.

%
% correspondence between E and E~.
%
\begin{pro}
Let $E \subset \BR^m$ and $\ti{E}$ be as above. If for 
$\veC{u}{l} \in (\BR^m)^l$ we have
\begin{equation}\label{eq:intersection}
  \mu(\ti{E} \cap T_{u_1}^{-1}\ti{E} \cap \ldots \cap T_{u_l}^{-1}\ti{E}) >0 , 
\end{equation}
then for all $\gd>0$,
\[
  E_{\gd} \cap (E_{\gd}-u_1) \cap \ldots \cap (E_{\gd}-u_l) 
  \ne \emptyset .
\]
\end{pro}

\begin{proof}
Define the function $g$ on $X^0$ by
\begin{equation*}
  g(\gp)=
   \begin{cases}
     \gd-f_0(\gp)& \te{if $f_0(\gp)< \gd$},   \\
     0           & \te{if $f_0(\gp) \ge \gd$}.
   \end{cases}
\end{equation*}
Since $g(\gp)$ is positive for $\gp \in \ti{E}$, equation 
(\ref{eq:intersection}) 
implies that
\[
  \int g(\gp)g(T_{u_1}\gp) \ldots g(T_{u_l}\gp) d\mu > 0 .
\]
In particular for some $\gp=T_w\gf$ the integrand is positive. As 
\[
  g(T_w \gf) > 0 \iff \gf(w) < \gd \iff  w \in E_{\gd} 
\]
we have 
\[
  w \in E_{\gd}, w+u_1 \in E_{\gd}, \ldots , w+u_l \in E_{\gd}.
\]
\end{proof}

We now forget the original set $E$, and the geometric problem takes the 
following dynamical form: 
\begin{thm}[Dynamical Version]\label{DynamicalVersion}
Let $(X,\caB,\mu,\BR^m)$ be a $\BR^m$ action, and let $T_u$ denote the action
of $u \in \BR^m$. Let $\veC{u}{k} \in (\BR^m)^k$, 
and let $A \subset X$, with $\mu(A) > 0$. There exists $t_0 \in \BR^{+}$ s.t. 
for all 
$t>t_0$, there exists a rotation $P \in SO(m)$ such that
\[
  \mu(A \cap T_{tP u_1}^{-1}A \cap \ldots \cap T_{tP u_{k}}^{-1}A) > 0 .
\] 
(Here $SO(m)$ is the special orthogonal group acting on $\BR^m$).
\end{thm}

%%%%%%%%%%%%%%%%%%%%%%%%%%% 
%%%%%%%%%%%%%%%%%%%%%%%%%%%
\section{Preliminaries.} %%
%%%%%%%%%%%%%%%%%%%%%%%%%%%
%%%%%%%%%%%%%%%%%%%%%%%%%%%

In the following section we give some measure theoretic and
ergodic theory preliminaries. The theorems are stated without proofs.
For the proofs see ~\cite{Fu1}, ~\cite{Pe}.

A \emph{measure preserving system} (m.p.s) is a system 
$X=(X^0,\caB,\mu,G)$ where $X^0$ is an arbitrary space, $\caB$ is a 
$\gs$-algebra of subsets of $X^0$, $\mu$ is a $\gs$-additive probability 
measure on the sets of $\caB$, and
$G$ is a locally compact group acting on $X^0$ by measure preserving 
transformations. We denote the action of the element $g \in G$ by $T_g$.
If the group $G =\BZ$, and $T$ is the generator of the $\BZ$ action, we
denote the system $(X^0,\caB, \mu,T)$.    
We say that the action of $G$ is \emph{ergodic}, if for any $A \in \caB$,
$T_g^{-1}A=A$ 
$\forall g \in G$, implies $\mu(A)=0$ or $\mu(A)=1$.
In this case we also say that $\mu$ is ergodic with respect to the 
action of $G$.
Each $T_g$ induces a natural operator on $L^2(X)$ by
$T_gf=f \circ T_g$, and the ergodicity of the action of $G$ is equivalent
to the assertion that there are no non-constant $G$-invariant functions.
We have:
\begin{thm}[Mean Ergodic Theorem]
Let $X=(X^0,\caB,\mu,T)$ be a m.p.s., \\
and $f \in L^2(X)$. Then
\[
  \Avr{N}{n} f \circ T^n \overset{L^2(X)}{\longrightarrow} \BP f,
\]
where $\BP f$ is the orthogonal projection of $f$ on the subspace of the 
$T$-invariant functions. 
\end{thm}

Let ${X}=(X^0,\caB,\mu,G)$ 
be a measure preserving system (m.p.s).
Let ${Y}=(Y^0,\caD,\nu,G)$ 
be a homomorphic image of $X$; i.e.,  we have a map  $\pi:X^0 \ra Y^0$
with $\pi^{-1}\caB_Y \subset \caB_X$, $\pi \mu_X =\mu_Y$ and $\pi$ commutes with the $G$ action.  
Then ${Y}$ is  a {\em factor} of ${X}$, ${X}$
is an {\em extension} of ${Y}$, and abusing the notation we write
$\pi:X \ra Y$ for the factor map. A factor of $X$ is determined by 
a $G$-invariant subalgebra of $L^{\infty}(X)$.
The map $\pi$ induces two natural maps
$\pi^*: L^2(Y) \ra  L^2(X)$ given by $\pi^*f=f\circ \pi$, and 
$\pi_*: L^2(X) \ra  L^2(\mu_{Y})$ given by $\pi_*f=E(f|\caB_{Y})$
(the orthogonal projection of $f$ on $\pi^*L^2(Y)$). The two measure 
preserving systems are \emph{equivalent} if the homomorphism 
of one to the other is invertible. We shall simplify the the notation writing $E(f|Y)$ for $E(f|B_Y)$.

A m.p.s. $X$ is \emph{regular} if $X^0$ is a compact
metric space, $\caB$ the Borel algebra of $X^0$, $\mu$ a measure on $\caB$.
A m.p.s. is \emph{separable} if $\caB$ is generated by a countable subset.
As every separable m.p.s. is equivalent to a regular m.p.s.,
we will confine our attention to regular m.p.s.\vspace{.10in}
 
% Disintegration of Measure
%
\subsection{Disintegration of Measure}
Let $(X^0,\caB,\mu)$ be a regular measure space, and let  
$\ga : (X^0,\caB,\mu)\ra (Y^0,\caD,\nu)$ be a homomorphism to 
another measure space (not necessarily regular). Suppose $\ga$ is induced 
by a   map $\gf: X^0 \ra Y^0$. In this case the measure $\mu$ has a 
{\em disintegration} in terms of fiber measures $\mu_y$, where $\mu_y$ is concentrated
on the fiber $\gf^{-1}(y)=X_y$. 
We denote by $\caM(X)$ the compact 
metric space of probability measures on $X^0$.
\begin{thm}\label{T:Disintegration}
There exists a measurable map from $Y^0$ to $\caM(X^0)$, $y \ra \mu_y$ 
which satisfies:
  \begin{enumerate}
    \item For every $f \in L^1(X^0,\caB,\mu)$, $f \in L^1(X^0,\caB,\mu_y)$ \\  
          for a.e. $y \in Y^0$, and $E(f|Y)(y) = \int f d\mu_y \quad$ for 
          a.e. $y \in Y^0$                                              
    \item $\int \{ \int f d\mu_y \} d \nu(y) = \int f d\mu \quad$
          for every $f \in L^1(X^0,\caB,\mu)$.
  \end{enumerate}
\end{thm}
The map $y \ra \mu_y$ is characterized by condition (1).
We shall write $\mu= \int \mu_y d\nu$ and refer to this as the 
disintegration of the measure $\mu$ with respect to $\caD$.  
 
If $(X^0,\caB,\mu,G)$  is a m.p.s., $\caD$ the algebra of all $G$-invariant sets, 
$\mu= \int \mu_x d\mu(x)$ the disintegration of $\mu$ with respect to $\caD$,
then $\mu_x$ is $G$-invariant and ergodic for a.e. $x$. 

%
%Definition of Universal Characteristic factors
%
\begin{dsc}{\bf Nilsystems and Characteristic Factors}
A {\em $k$-step nilflow} is a system $X=(N/\gC,\caB,m,G)$ where
$N$ is a $k$-step nilpotent Lie group, $\gC$ a cocompact lattice,
$\caB$ the (completed) Borel algebra, $m$ the Haar measure, and the action 
of $G$ is by 
translation by elements of $N$: $T_gn\gC=a_gn\gC$ where $g \ra a_g$ is a 
homomorphism of $G$ to $N$.
We will sometimes denote this system by  $(N/\gC,G)$, or $(N/\gC, a)$ if 
$G=\BZ$ and $1 \ra a$. If $G$ is connected
and $(N/\gC,G)$ is an ergodic nilflow, then we may assume that $N$ is connected
so that $X^0 = N/\gC$ is connected and is a homogeneous space of the identity 
component of $N$.
A {\em $k$-step pro-nilflow} is an inverse limit of $k$-step nilflows.
\end{dsc}

\begin{thm}[Cf.\cite{Pa1}] \label{thm:uniform}
Let $X=(N/\gC,a)$ be an ergodic nilflow,  then 
$X$ is uniquely ergodic. Let $f$ be a continuous 
function on $N/\gC$. Then the averages $\Avr{N}{n} f(a^nx)$ converge 
uniformly to $\int f(x) dm$.
\end{thm}

\begin{dsc}\label{dsc:sasha}
Let $N$ be a connected simply connected nilpotent Lie group, $\gC$ a cocompact
lattice 
in $N$, and $X^0=N/\gC$. Let $\pi:N \ra X^0$ be the natural projection, 
and let $M$ be a closed connected subgroup of $N$ such that $\pi(M)$ is a 
closed submanifold of $X^0$. Let $G=\BR^k$ and let $\gf:G \ra N$ be 
a homomorphism. For $x \in X^0$ let $O(x)=\overline{Gx}$, and for $x \in X^0$, 
$g \in G$ let 
$O_g(x)=\overline{\{\gf(ng)x\}}_{n \in \BZ}$; these are subnilmanifolds of 
$X^0$ (see for example \cite{Le}) .

Introducing Malcev coordinates on $N$ and $M$ (\cite{Ma}) 
we can identify these groups topologically with, say, $\BR^l$ and $\BR^m$, 
$l \ge m$. 
Call a proper subspace of $\BR^d$ 
{\em countably linear} if it is contained in a countable union of proper 
linear subspaces. Call a subset of $\BR^d$ {\em polynomial} if it is 
the set of zeroes of some nonzero polynomial in $\BR^d$ (i.e. an algebraic 
variety of co-dimension $1$), and {\em countably 
polynomial} if it is contained in a countable union of proper polynomial 
subsets. The following proposition is due to Sasha Leibman: 

\begin{pro}\label{pro:sasha} There exists a connected subnilmanifold $V$ of 
$X^0$ such that
$O(\pi(a)) \subseteq aV$ for all $a \in M$, and there exists a countably linear
set $B \subset G$ such that for every $g \in G \setminus B$ there is a countably
polynomial set $A_g \subset M$ such that $O_g(\pi(a))=aV$ for all 
$a \in M \setminus A_g$.
\end{pro} 

\begin{proof} Define a mapping $\eta:G \times M \ra N$ by
$\eta(g,a)=a^{-1}\gf(g)a$. In Malcev coordinates on $M$ and $N$,
$\eta$ is a polynomial mapping $\BR^{k+m} \ra \BR^l$ (see \cite{Ma}). Moreover 
for each $a \in M$, $\eta(\cdot,a)$ is a homomorphism $G \ra N$. 
Let $H$ the closure of the subgroup generated by $\eta(G \times M)$.
Let $V$ be the closure of $\pi(H)$ in $N/\gC$. Then  $V$ is a subnilmanifold 
$V=\pi(K)$ for some closed 
subgroup $K$ of $N$ (\cite{Sh}) ($\pi(H)$ itself is not necessarily closed).
We then have $a^{-1}\gf(g)\pi(a)=\pi(a^{-1}\gf(g)a) \in V$, thus 
$\gf(g)\pi(a) \in aV$ for any $a \in M$, and  $g \in G$. So, 
$O(\pi(a))\subseteq aV$ for all $a \in M$.

Let $\ti{L}$ be the set $\{ l \in N: lV=V\}$. Then $\ti{L}$ is a group,
$\eta(G \times M)$ and $K$ are subsets of $\ti{L}$, and $V=\pi(\ti{L})$. 
Let $L$ be the identity component of $\ti{L}$, then 
$\eta(G \times M) \subseteq L$, and therefore $H  \subseteq L$.
$V$ is connected, and therefore a homogeneous subspace of $L$; $V=L/L\cap\gC$.
Let $W$ be the maximal torus factor of $V$, 
$W=L/([L,L](L\cap\gC))$, and let $p:V \ra W$ be the natural projection. 
Let $\hat{W}$ be the group of characters of $W$, and let $\gx \in \hat{W}$. 
The character $\gx$ can be lifted to a homomorphism $\gz_{\gx}:L \ra \BR$. 
For each $\gx \in \hat{W}$, let $\gp_{\gx}:=\gz_{\gx} \circ \eta$. 
Then $\gp_{\gx}$ are polynomials
on $G \times M$, which for each $a \in M$ are linear with respect to $G$. Moreover, each $\gp_{\gx}$ is a nonzero polynomial; otherwise 
$\eta(G \times M)$ would be contained in the kernel of the corresponding 
homomorphism  $\gx \circ p \circ \pi : N \ra S^1$. This is a closed 
subgroup of $N$ containing $\eta(G \times M)$, thus contains the subgroup $H$.
Therefore  $\gx \circ p \circ \pi (H)=1$, but this implies that 
$\gx \circ p (V) =1$, i.e. $\gx$ is the trivial character.

Let $C_{\gx} \subset G \times M$ be the set of zeros of $\gp_{\gx}$, and let 
$C=\bigcup_{\gx \in \hat{W}} C_{\gx}$. Then
$C$ is a countably polynomial subset of $G \times M$.

For any $(g,a) \notin C$ one has 
$\gx \circ p(a^{-1}\gf(g)\pi(a))=\gx \circ p (\pi(a^{-1}\gf(g)a))\ne 0$ 
for all $\gx \in \hat{W}$, so the projection of 
$a^{-1}\gf(g)\pi(a) \in V$ to $W$ is not 
contained in any proper 
subtorus of $W$. Consider the following $\BZ$ action on $V$: for  $n \in \BZ$, 
$v \ra a^{-1}\gf(ng)av$. Since the projection of 
$a^{-1}O_g(\pi(a)) \subseteq V$
to $W$ is a closed subgroup of $W$, i.e. a subtorus of $W$, 
it is equal to $W$.   By Parry (\cite{Pa1}) this implies  the $\BZ$ action is 
minimal and therefore
$a^{-1}O_g(\pi(a))=V$, and so $O_g(\pi(a))=aV$.

Now let
\[
B=\{ g \in G:\{g\} \times M \subseteq C\}, \  \
M_{\gx}(g)= \{a \in M: \gp_{\gx}(g,a)=0 \}.
\]
If $\{g\} \times M \subseteq C$, then 
$M= \bigcup_{\gx \in \hat{W}}  M_{\gx}(g)$. As $M$ is connected, 
if $M_{\gx}(g)$ has non-empty interior, then $M_{\gx}(g)=M$. 
By the Baire category theorem $M_{\gx}$ 
is non-empty for some $\gx \in \hat{W}$. Therefore 
\[
B =\bigcup_{\gx\in \hat{W}} \{ g \in G: \gp_{\gx}(g,a)=0 \ 
\te{for all} \ a \in M\}.
\]
Then $B$ is a countably linear subset of $G$, and for each 
$g \in G \setminus B$, 
\[
A_g=C \cap (\{g\} \times M) = \bigcup_{\gx \in \hat{W}} 
  \{a \in M: \gp_{\gx}(g,a)=0\}
\] 
is a countably polynomial subset of $M$.
\end{proof}
\end{dsc}

\begin{thm}\label{thm:abelian_structure}  Let $X=(X^0,\caB,\mu,\BR^m)$ be an ergodic $\BR^m$ action. We can associate with
$X$ an inverse sequence of factors  
$\ldots \ra Y_k(X) \ra Y_{k-1}(X) \ra \ldots \ra Y_1(X)$, where $Y_k(X)$ is 
a $k-1$-step
pro-nilflow such that the following holds: 
If $u_1,\ldots,u_k \in \BR^m$ are such that the actions
of $T_{u_i}$ and $T_{u_i-u_j}$ for $i \ne j$ are ergodic, 
then for any bounded measurable functions $f_1,\ldots,f_k$ the limits
in $L^2(X)$ 
\begin{equation}\label{eq:FK1}
\lim_{n \ra \infty}  \Avr{N}{n} \p_{j=1}^k  T_{nu_j} f_j(x), \ 
\lim_{n \ra \infty} \Avr{N}{n} \pi^*\p_{j=1}^k T_{nu_j} E(f_j|Y_k)(x) 
\end{equation}
exist and are equal. The factor $Y_k(X)$ is called the {\em $k$-universal 
characteristic factor ($k$-u.c.f) of $X$} .
Let
\[ \begin{split}
   \tau_{\vec{u}}(T):=& T_{u_1} \times \ldots \times T_{u_l},\\
   \end{split}
\] 
let $\triangle_{k}(\gm)$ be the diagonal measure on $X^k$, then 
\[        \bar{\triangle}_{\vec{u}}(\mu):= 
        \lim_{N \ra \infty} \Avr{N}{n} \tau_{\vec{u}}^n \triangle_{k}(\mu)
\]
is well defined. If $F$ is a function
invariant under $\tau_{\vec{u}}$ with respect to the measure 
$\bar{\triangle}_{\vec{u}}(\mu)$
and if $E(f_j|\caB_{Y_k})=0$ for some $1\le j \le k$, then 
\[
\int f_1(x_1)\ldots f_k(x_k)F(x_1,\ldots,x_k)d\bar{\triangle}_{\vec{u}}=0.
\]  
\end{thm}

The factors $Y_k(X)$ were constructed for an ergodic m.p.s
$X=(X^0,\caB,\mu,T)$ by Host and Kra \cite{HKr} and independently by Ziegler
\cite{Z}. Frantzikinakis and Kra \cite{FrKr} showed that  if 
$X_i=(X^0,\caB,\mu,T_i)$ are ergodic measure 
preserving systems on the same space $X^0$, where $T_i$ commute, then 
$\caB_{Y_k(X_i)}=\caB_{Y_k(X_j)}=\caB_{Y_k}$ for any $i,j$, and if the action of 
$T_i^{-1}T_j$ is ergodic for all $i \ne j$, then 
equation (\ref{eq:FK1}) holds (replacing $T_{nu_i}$ with 
$T_{ni}=T^n_i$). Thus if we have a
$\BR^m$ action then the systems $X_u=(X^0,\caB,\mu,T_u)$, $u \in \BR^m$
for which the action of $T_u$ is ergodic, share the same sequence
of factors $Y_k(X)=Y_k(X_u)$. The fact that $Y_k(X)$ is a factor of the 
$\BR^m$ action follows from \cite{Z} corollary $2.4$. We will show that the
action of $\BR^m$ on $Y_k(X)$ preserves the pro-nil structure: 

\begin{dfn}
Let $Y=(Y^0,\caB_Y,\mu_Y,\BR^m)$ be a $j$-step pronilflow; 
$Y=\lim_{\leftarrow} N_i/\gC_i$. We say that the action of $\BR^m$ on  $Y$ 
preserves the  pro-nil structure 
if the action of $\BR^m$ on $Y$ induces a $\BR^m$ action on
$N_i/\gC_i$ by group rotations. 
\end{dfn}

%%%%%%%%%%%%%%%%%%%%%%%%%%%%%%%%%%
\begin{pro} Let $Y=(Y^0,\caB,\mu,T)$ be a $j$-step ergodic pronilflow;
$Y=\lim_{\leftarrow} N_i/\gC_i$. Let $\{T_c\}_{c \in \BR^m}$ be a
$\BR^m$ action on  $(Y^0,\caB,\mu)$ that 
commutes with the action of $T$. Then the action of $\BR^m$ on $Y$ 
preserves the pro-nil structure.
\end{pro}
%%%%%%%%%%%%%%%%%%%%%%%%%%%%%%%%%

\begin{proof}
For $j=1$, $Y$ is a Kronecker action, and any factor of $Y$ is a 
Kronecker action. Thus it is enough to check that eigenfunctions of the 
$T$ action are also eigenfunctions of the $\BR^m$ action.
If $\psi$ is an eigenfunction, $T\psi(y)=\gl\psi(y)$, then as $T$ and $T_c$
commute $\psi(TT_cy)=\gl \psi(T_cy)$. Combining the two equations we get
\[
T\left( \frac{\psi(T_cy)}{\psi(y)}\right)=1.
\]
By ergodicity of $T$ we get $\psi(T_cy)=\gd_c \psi(y)$.

We proceed by induction on $j$. Let $Y$ be a $j$-step ergodic pronilflow; 
$Y= \lim_{\leftarrow} M_i/\gL_i$. We first show that the $\BR^m$ action
on Y induces a  $\BR^m$ action on $M_i/\gL_i$. Let $\pi:Y \ra M_i/\gL_i$ be 
the projection. Let $p: Y \ra Y_j(M_i/\gL_i)$ be the projection onto the 
$j$ u.c.f of $M_i/\gL_i$. $Y_j(M_i/\gL_i)$ is a $j-1$-step nilflow, we denote
it $N_i/\gC_i$. 
The space $L^2(M_i/\gL_i)\circ \pi \subset L^2(Y)$ is 
spanned by functions $f$ satisfying the following condition:
\[
Tf(y)=g(y)f(y)
\]
where $g=\ti{g} \circ p$ with  $\ti{g}$ of type $j$ 
(see \cite{Z} theorem $6.1$).
As $T,T_c$ commute for any $c \in \BR^m$ 
\[
TT_cf(y)=T_cTf(y)=T_cg(y)T_cf(y).
\]
Thus
\[
T\left( \frac{f(T_cy)}{f(y)}\right)=\frac{T_cg(y)}{g(y)}\frac{f(T_cy)}{f(y)}.
\]
By the induction hypothesis the action of $\BR^m$ on $Y$ induces an action
on $Y_j(M_i/\gL_i)=N_i/\gC_i$, and this action is given by rotation
by an element $a_i(c) \in N_i$.
By  proposition $6.37$ in \cite{Z}, as $\BR^m$ commutes with the 
action of $T$ on $N_i/\gC_i$ given by rotation by $a \in N_i$, 
there exists a family of measurable functions  
$\{ f_{c}:N_i/\gC_i \ra S^1 \}_{c \in \BR^m}$ and a family of constants
$\{ \gl_c \}_{c \in \BR^m}$ such that 
 \[
\frac{T_c \ti{g}(p(y))}{\ti{g}( p(y))}
=\gl_{c}\frac{Tf_{c}(p(y))}{f_{c}(p(y))}.
\]
We get
\[
T\left( \frac{f(T_cy)}{f(y)f_c(p(y))}\right)=\gl_c
\frac{f(T_cy)}{f(y)f_c((p(y))}.
\]
This implies that $\gl_c$ is an eigenvalue of $T$, but as it is 
multiplicative in $c \in \BR^m$, $\gl_c \equiv 1$.
Therefore by ergodicity of the $T$ action
\[
\frac{f(T_cy)}{f(y)f_c(p(y))}=\gd_c'
\]
or
\[
f(T_cy)=\gd_c f(y) f_c(p(y)) \in L^2(M_i/\gL_i)\circ \pi.
\]
This shows that the $\BR^m$ action on $Y$ induces an $\BR^m$ action
on $M_i/\gL_i$. The fact that this action is given by group rotations
was shown by Parry \cite{Pa2} in the case where $M_i$ is connected.
Alternatively, $M_i/\gL_i$ can be presented as a torus extension of 
$Y_j(M_i/\gL_i)=N_i/\gC_i$ with
$g:N_i/\gC_i \ra \BT^n$ a cocycle of type $j$. Without loss of generality
we can assume $n=1$. Now the tuples $(a,g)$, $(a_i(c),f_c)$ belong to 
the group 
$\caG$ defined in \cite{Z}  proposition $6.37$ and this group acts 
transitively and effectively on $M_i/\gL_i$.
\end{proof}

\begin{pro}[\cite{PS}]\label{P:ergodic}
If $(X,\caB,\mu,\BR)$ is an ergodic action of $\BR$, then but for a countable
set of $u \in \BR$, $T_u$ acts ergodically.
If $(X,\caB,\mu,\BR^m)$ is an ergodic action of $\BR^m$, then but for
a countable set of $l-1$ dimensional hyperplanes,
all $l-1$ dimensional hyperplanes through the origin act ergodically.
\end{pro}

The following is a version of the van der Corput Lemma 
(see \cite{FuKaW}).
%
% Lemma hilbertian
%
\begin{lma}\label{L:hilbert}
Let $H$ be a Hilbert space, $\xi \in \Xi$ some index set,  
and let $u_n(\xi) \in H$  for $n \in \BN$ be uniformly bounded in $n,\xi$. 
Assume that for each $r$ the limit
\[
  \gc_r(\xi)= \Lim{N} \Avr{N}{n} \ip{u_n(\xi)}{u_{n+r}(\xi)}
\]
exists uniformly and
\begin{equation}
  \Lim{R} \Avr{R}{r} \gc_r(\xi)=0
\end{equation}
uniformly. Then
\[
  \Avr{N}{n} u_n(\xi) \overset{H}{\lra} 0
\]
uniformly in $\xi$.
\end{lma}

\begin{dsc}{\bf Multidimensional Szemer\'edi.}
The following generalization of Szemer\'edi's theorem was proved 
by Furstenberg and Katznelson  ~\cite{FuKa}:
\begin{thm}\label{thm:FuK}
Let $X=(X^0,\caB,\gm,\BZ^k)$ be a m.p.s., and let $T_1,\ldots,T_k$ be the 
generators of the $\BZ^k$ action. 
Let $f\ge 0$ be a bounded measurable function on $X$ with $\int f d\mu >0$. Then 
\[
\liminf_{N \ra \infty}  \Avr{N}{n} \int f(x)T^n_1f(x) \ldots T^n_kf(x)
d\mu(x) >0.
\]
\end{thm}
\end{dsc}

%%%%%%%%%%%%%%%%%%%%%%%%%%%%%%%%%%%%%%%
%%%%%%%%%%%%%%%%%%%%%%%%%%%%%%%%%%%%%%%
\section{The Main Theorem}           %%
%%%%%%%%%%%%%%%%%%%%%%%%%%%%%%%%%%%%%%%
%%%%%%%%%%%%%%%%%%%%%%%%%%%%%%%%%%%%%%%

Denote  $M_m(\BR)$ the $m \times m$ matrices over $\BR$, and $SO(m)$
the special orthogonal group. Recall that if $(N/\gC,G)$ is a nilflow
the action of $T_g$ for $g \in G$ is given by 
$T_gn\gC=a_gn\gC$ where $a_g \in N$.

\begin{lma} Let $(N/\gC,\BR^m)$ be an ergodic 
measure preserving  action of $\BR^m$ on a 
nilmanifold $N/\gC$, where $N$ is connected. 
Let $f_j$ be continuous functions on $N/\gC$. 
Let $(u_1,\ldots, u_l) \in (\BR^m)^l$. Then there exists a countably linear 
set (see \ref{dsc:sasha}) $\caS \subset M_m(\BR)$ such that for any $F \in M_m(\BR) \setminus \caS$ 
the function 
\begin{equation}\label{eq:independence}
g_{F,L}(x):=\lim_{N \ra \infty} \Avr{N}{n} \p_{j=1}^l T_{(nF+L)u_j} f_j(x)
\end{equation}
is independent of  $L \in M_m(\BR)$ for a.e. $x \in N/\gC$.
Furthermore for any such $F$ the convergence is uniform in $L$.
\end{lma}

\begin{proof} Let $M$ be the diagonal of $N^l$ and let $G=M_m(\BR) =\BR^{m^2}$ 
(thought of as an additive group). Let $\gf: M_m(\BR) \ra N^l$ 
be given by 
\[
\gf(F)=(a_{Fu_1},\ldots,a_{Fu_l}).
\]
By proposition 
\ref{pro:sasha}, there exists a submanifold $V$ of $(N/\gC)^l$, and there 
exists a countably linear set $\caS \subset M_m(\BR) $  such that for every 
$F \in M_m(\BR) \setminus \caS$ there is a countably polynomial set
$A_F \subset M$ 
such that for $(a,\ldots,a) \notin A_F$, 
\[ 
\overline{\{\gf(nF)\pi(a,\ldots,a)\}}_{n \in \BZ}=(a,\ldots,a)V,
\]
and
\[
\overline{G\pi (a,\ldots,a)} \subseteq (a,\ldots,a)V \ \ \
(\te{therefore} = (a,\ldots,a)V).
\] 
For any $F \in  M_m(\BR) \setminus \caS$,  and $(a,\ldots,a) \notin A_F$
we have  
\[
T_{Lu_1}\times \ldots \times T_{Lu_l} \pi(a,\ldots,a) \in 
(a,\ldots,a)V.
\]
The action of $\gf(F)$ on $(a,\ldots,a)V$ is ergodic, and by 
theorem \ref{thm:uniform} it is uniquely ergodic. The point 
$(T_{Lu_1}a\gC,\ldots,T_{Lu_l}a\gC) \in (a,\ldots,a)V$.   
By theorem \ref{thm:uniform} the convergence in equation
(\ref{eq:independence}) is uniform in $L$, and $g_{F,L}(a\gC)$ is 
independent of  $L$.
\end{proof}

\begin{cor}\label{cor:same_limit} Let $Y=(Y^0,\caB,\mu,\BR^m)$ 
be an ergodic pro-nilflow.
Let $f_j$ be bounded measurable functions on $Y^0$. 
Let $(u_1,\ldots, u_l) \in (\BR^m)^l$. 
Then there exists a countably linear 
set $\caS \subset M_m(\BR)$ such that for any $F \in M_m(\BR) \setminus \caS$,
and all $L \in M_m(\BR)$
the function
\[
g_{F,L}(y):=\lim_{N \ra \infty} \Avr{N}{n} \p_{j=1}^l T_{(Fn+L)u_j} f_j(y)
\]
where the limit is in $L^2(Y)$, is a constant function of $L \in M_m(\BR)$ 
and the convergence 
is uniform in $L$.
\end{cor}

\begin{proof} If $Y^0=\lim_{\leftarrow} N_j/\gC_j$, the continuous functions on 
$N_j/\gC_j$ lifted to $Y^0$, for all $j$, are dense in $C(Y^0)$.
\end{proof}

The next proposition will enable us to evaluate averages of functions
on $X$ by evaluating the averages of the projections of the functions
on the factor $Y_k(X)$ described in \ref{thm:abelian_structure}. 
 
\begin{pro}\label{P:AbelianReduction}
Let $X=(X^0,\caB,\gm,\BR^m)$ be an ergodic action of $\BR^m$, 
and let  $\veC{u}{k} \in (\BR^m)^k$.
Let
$Y_k$ be the factor described in theorem \ref{thm:abelian_structure}, 
and let $\pi:X \ra Y_k$ be the factor map. 
Let $f_1 \ldots f_k$ be bounded measurable functions on $X$. 
Then there exists a countably linear subset $\caS \subset M_m(\BR)$
such that for any $M \in M_m(\BR) \setminus \caS$, satisfying
$T_{Mu_i}$, $T_{M(u_i-u_j)}$  
for $i,j=1,\ldots,k$, $i \ne j$ are ergodic, and for all $P \in M_m(\BR)$
we have
\[
   \Avr{N}{n} \p_{j=1}^k T_{(nM+P)u_j} f_j(x) - 
   \Avr{N}{n} \p_{j=1}^k T_{(nM+P)u_j} \pi^*E(f_j|Y_k)(x)  
   \overset{L^2(X)}{\lra} 0 
\]
uniformly in $P$.
\end{pro}

\begin{proof}
We prove this inductively. For $k=1$, let $u \in \BR^k$, $u \ne 0$.
If $T_{Mu}$ is ergodic then
\[
   \Avr{N}{n} T_{nMu+Pu} f(x)= T_{Pu}(\Avr{N}{n} 
   T_{nMu}f(x))\lra \int f(x) d\mu
\]
uniformly in $P$, by the Mean Ergodic Theorem. 
Assume the statement holds for $k$: i.e., for $M$ outside a countably linear set
satisfying  $T_{Mu_i}$, $T_{M(u_i-u_j)}$ for $i,j=1,\ldots,k$,  
$i \ne j$ are ergodic, and all $P \in  M_m(\BR) $ we have 
\[
   \Avr{N}{n} \p_{j=1}^k T_{(nM+P)u_j} f_j(x) - 
   \Avr{N}{n} \p_{j=1}^k T_{(nM+P)u_j} \pi^*E(f_j|Y_k)(x)  
   \overset{L^2(X)}{\lra} 0 
\]
uniformly in $P$.
We show this for $k+1$.
Let $\caS \subset M_m(\BR)$ be the set from corollary \ref{cor:same_limit}
corresponding to $Y_k$ and $u_1, \ldots, u_{k+1}$. For 
$M\in   M_m(\BR) \setminus \caS$ the $L^2$ limit 
\[
\lim_{N \ra \infty} \Avr{N}{n} \p_{j=1}^{k+1} T_{(nM+P)u_j} E(f_j|Y_k)(y)
\] 
is independent of $P$, and the convergence to the limit in uniform in $P$.
Let $M \in  M_m(\BR)\setminus \caS$ satisfy $T_{Mu_i}$,  $T_{M(u_i-u_j)}$
are ergodic for $i,j=1,\ldots,k+1$, $i \ne j$.
It is enough to show that if 
for some $1 \le j \le k+1$,  $E(f_j|Y_{k+1})=0$ then 
\[
\lim_{N \ra \infty} \Avr{N}{n} \p_{j=1}^{k+1} T_{(nM+P)u_j} f_j(x)=0
\]
uniformly in $P$.
We use the Van der Corput Lemma (lemma \ref{L:hilbert}).
Let $v_n(M,P):= \p_{j=1}^{k+1} T_{nMu_j+Pu_j} f_j(x)$.
Then 
\[ 
\ip{v_n(M,P)}{v_{n+r}(M,P)}=\int \p_{j=1}^{k+1} T_{nMu_j+Pu_j} f_j(x)T_{(n+r)Mu_j+Pu_j} \bar{f}_j(x) \ d\mu,
\]
and 
\begin{equation*}
\begin{split}
\gc_r &(M,P) =  \lim_{N \ra \infty}  \Avr{N}{n}\ip{v_n(M,P)}{v_{n+r}(P,M)} \\
                  = & \lim_{N \ra \infty} \int f(x)T_{rMu_1}\bar{f}(x) 
                      \Avr{N}{n}\p_{j=2}^{k+1}T_{nM(u_j-u_1)+P(u_j-u_1)} 
                      (f_j(x)T_{rMu_j}\bar{f}_j(x))\ d\mu. \\
\end{split}
\end{equation*}
By the induction hypothesis this limit is equal (uniformly in  $P$) 
to the following limit
\begin{equation}\label{eq:Y_K}
 \lim_{N \ra \infty}\int f(x)T_{rMu_1}\bar{f}(x) 
                      \Avr{N}{n}\p_{j=2}^{k+1}T_{nM(u_j-u_1)+Pu_j-Pu_1} 
                      \pi^*E(f_jT_{rMu_j}\bar{f}_j|Y_{k})(x)\ d\mu,
\end{equation}
which equals
\begin{equation}\label{eq:Y_K_1}
 \lim_{N \ra \infty} \Avr{N}{n}\p_{j=1}^{k+1}
                    \int T_{nMu_j+Pu_j}
                    \pi^*E(f_jT_{rMu_j}\bar{f}_j|Y_k)(x)\ d\mu.
 \end{equation}
The limit in equation (\ref{eq:Y_K_1})  is a limit on a 
$k-1$ step pronilflow, as $M \in M_m(\BR) \setminus \caS$  it is the same for 
all $P$, and the convergence is uniform in $P$.
By \ref{thm:abelian_structure}, the limit in equation (\ref{eq:Y_K_1}) is equal to  
\[
     \int \p_{j=1}^{k+1} T_{Pu_j}f_j(x_j)T_{rMu_j+Pu_j}\bar{f}_j(x_j)
                      \ d\bt_{M\vec{u}}(\mu)(x_1,\ldots,x_{k+1}),
\]
where $\bt_{M\vec{u}}(\mu)$ is a measure on $X^{k+1}$.
Now 
\[
\lim_{R \ra \infty}\Avr{R}{r} \p_{j=1}^{k+1} T_{rMu_j+Pu_j}\bar{f}_j(x_j)
\]
converges uniformly in $P$ to a function $F$ in 
$L^2(\bt_{M\vec{u}}(\mu))$ which is invariant under 
$T_{Mu_1} \times \ldots \times T_{Mu_{k+1}}$ 
(by the Mean Ergodic Theorem, as in the case $k=1$). 
Finally by  \ref{thm:abelian_structure}, if $E(f_j|Y_{k+1})=0$ for some 
$1 \le j\le k+1$, then 
\[
\lim_{R \ra \infty}\Avr{R}{r} \gc_r(M,P) =0.
\]
(uniformly in $P$).
\end{proof}

\begin{rmr}\label{rmr:complex}
Corollary \ref{cor:same_limit} and proposition \ref{P:AbelianReduction} remain 
valid if we replace $M_m(\BR)$ by a linear subspace of $M_m(\BR)$. We apply 
this for the case $m=2$, replacing $M_2(\BR)$ by the embedding 
$\BC \hra M_2(\BR)$. Then, thinking of $u_1, \ldots, u_k$ as points in $\BC$
we can replace the matrices $F,L \in M_2(\BR)$ with  
$c,d \in \BC$ where $c$ is outside countably many lines in $\BC$.
\end{rmr}

\begin{lma}\label{Bantisymetric}
For each $r=1,\ldots \infty$, let $\{s^r_{l}\}_{l=1}^m \subset \BR^m$, 
such that for each $r$, $s_l^r\neq \vec{0}$ for some $1\le l \le m$.
Let $\{e_l\}_{l=1}^m$ be the standard basis for $\BR^m$.
There exists an antisymmetric matrix $B \subset M_m(\BR)$ s.t. 
$Bu_i,B(u_i-u_j)\neq \vec{0}$ for $1\le i,j \le k$, $i \neq j$, and 
\[
 \forall r : f_{r,B}(M) \df \Sum{l}{m} \ip{s^r_{l}}{MBe_l} \not \equiv 0. 
\]
\end{lma}
\begin{proof}
Let $\caB$ be the subspace of antisymmetric matrices. Since $f_{r,B}(M)$
is linear in $M$, we have $f_{r,B}(M) \equiv 0 \iff B$ satisfies 
the $m^2$ linear equations 
given by the standard basis for $\BR^{m^2}$. Hence for each $r$, 
the 'bad' $B$ form a linear subspace of $\caB$. Since we have only
a countable number of inequalities, it suffices to show that this
linear subspace is a proper subspace of $\caB$. So without loss of
generality, we have only one inequality.    
Assume 
\[ 
  \forall B \in \caB :  \Sum{l}{m} \ip{Ms_l}{Be_l} \equiv 0 \qquad  
\]  
Without loss off generality $s_{11} \neq 0$.
Let $E_{21}$ be an $m \times m$  matrix with $1$ at the index $21$, and $0$ 
elsewhere.
Then
\[
\begin{split}
\Sum{l}{m} \ip{E_{21}s_l}{Be_l}&=s_{11} b_{21} +s_{21} b_{22}+\ldots+s_{m1} b_{2m}\\
                               &=-s_{11} b_{12} +s_{31} b_{23}+\ldots+s_{m1} b_{2m}=0
\end{split}
\]
As $s_{11}\neq 0$ this is a non trivial linear condition on antisymmetric matrices.
Finally, the conditions $Bu_i=\vec{0}$ or $,B(u_i-u_j)=\vec{0}$ are non trivial 
linear conditions on antisymmetric matrices.
\end{proof}

\begin{lma}\label{lma:antisymmetric}
Let $\caS$ be a countably linear set in $M_m(\BR)$, $m \ge 3$.
Let $u_1,\ldots,u_k \in \BR^m$. 
There exist matrices
$M \in M_m(\BR) \setminus \caS$, and $P \in SO(m)$ such that 
$M^tP$ is antisymmetric, and $T_{Mu_i}$, $T_{M(u_i-u_j)}$ 
for $i,j=1,\ldots,k$, $i \ne j$ are ergodic.
\end{lma}

\begin{proof}
The set $\caS$ is countably linear therefore it is a countable union of sets of 
the form
\[
\caS_r=\{ N \in M_m(\BR): \sum_{l=1}^m \ip{s^r_l}{Ne_l}=c_r\},
\]
Where $e_l$ is the standard basis for $\BR^m$, $s^r_l \in \BR^m$, $c_r \in \BR$..
By lemma \ref{Bantisymetric} there exists an antisymmetric matrix $B$, such that
$Bu_i,B(u_i-u_j)\neq \vec{0}$ for $1\le i,j \le k$, $i \neq j$,  and for all $r$
\[
f_r(M)=\sum_{l=1}^m \ip{s^r_l}{MBe_l} \not \equiv 0.
\]
For each $r$, the set of $M$ with $f_r(M)=c_r$ is a hyperplane in $M_m(\BR)$.
This subspace intersects $SO(m)$ in a proper algebraic subvariety of $SO(m)$.
Therefore for a.e. $P \in SO(m)$ (with respect to the Haar measure on $SO(m)$), 
$M=PB$  will avoid the bad set  $\tilde{\caS}$. Clearly, if $M=PB$ avoids $\caS$
then $tM=tPB$ avoids $\caS$ for any $t>0$. 
By proposition \ref{P:ergodic}, for a.e. $P \in SO(m)$ and a.e. $t \in \BR$
$T_{tPBu_i}$, $T_{tPB(u_i-u_j)}$ act ergodically.
\end{proof}

\begin{dsc}{\em Proof of theorem \ref{DynamicalVersion}}.
Without loss of generality, we may assume by disintegration of $\mu$, 
that the action of $\BR^m$ is ergodic.
Let $f=1_A$ be the characteristic function of the set $A$, and $\mu(A)=\gl$.
Let $Y_k$ be the factor described in theorem \ref{thm:abelian_structure}, and let  
$E(f|Y_k)$ be the projection of $f$ on $L^2(Y_k)$. 

We first prove the theorem for $m>2$.
By corollary \ref{cor:same_limit}, proposition \ref{P:AbelianReduction},
and lemma \ref{lma:antisymmetric}, there exist matrices  $M\in M_m(\BR)$, $P \in SO(m)$ such that 
$M^{t}P$ is antisymmetric, and for all $t \in \BR$  we have
\[
\lim_{N \ra \infty} \Avr{N}{n} \p_{j=1}^k T_{nMu_j+tPu_j} E(f|Y_k)(y)=g(y)
\]
in $L^2(Y_k)$, and 
\[
\Avr{N}{n} \p_{j=1}^k T_{nMu_j+tPu_j} f(x)-
\Avr{N}{n} \p_{j=1}^k T_{nMu_j+tPu_j} \pi^*E(f|Y_k)(x)\ra 0
\]
in $L^2(X)$, where the convergence is uniform in $t$.  Then 
\[
\lim_{N \ra \infty} \Avr{N}{n} \p_{j=1}^k T_{nMu_j+tPu_j} f(x)= 
\lim_{N \ra \infty}\Avr{N}{n} \p_{j=1}^k T_{nMu_j}f(x) = \pi^*g(x),
\]
and the convergence is uniform in $t$.
By theorem \ref{thm:FuK},
\[
\int f(x) \pi^*g(x) d\mu_X > C>0.
\]
Uniform convergence implies that there exists $N_0$, such that for all $t$
\[
\left|\Avr{N_0}{n} \int f(x) \p_{j=1}^k T_{nMu_j+tPu_j} f(x) d\mu_X
- \int f(x) \pi^*g(x) d\mu_X \right|< \frac{C}{2}.
\]

Therefore for $N_0$, and for all $t \in \BR$
\[
\Avr{N_0}{n} \int f(x) \p_{j=1}^k T_{nMu_j+tPu_j} f(x) d\mu_X >\frac{C}{2}.
\]
This implies that for all  $t \in \BR$ there exists $n\le N_0$ 
with  
\begin{equation*}
\begin{split}
\mu( A\cap & T_{(nM+tP)u_1}A\cap \ldots \cap T_{(nM+tP)u_k}A) = \\
&\int f(x)\p_{j=1}^k T_{(nM+tP)u_j} f(x) d\mu_X 
>\frac{C}{2}.
\end{split}
\end{equation*}
Now the $T_u$ satisfy the following continuity condition:
\begin{equation} \label{eq:cont}
  \forall \gve \exists \gd : \|u-u'\| \le \gd \Rightarrow
  | \mu(A \cap T_uA) - \mu(A \cap T_{u'}A) | \le \gve.
\end{equation}  
As $M^tP$ is antisymmetric, $M \in T_P(SO(m))$ - the tangent space of $SO(m)$
at $P$. Thus 
\[
P':=Pexp(\epsilon nP^{-1}M)=P(I+\epsilon nP^{-1}M+o(\epsilon))=
   P+\epsilon n M + o(\epsilon)
\]
belongs to $SO(m)$.
But
\[
(\frac{1}{\epsilon}P+nM)- \frac{1}{\epsilon}P'=o(1),
\]
and if $t=\frac{1}{\epsilon}$ is large enough, then by equation 
(\ref{eq:cont}) 
\[
\mu(A \cap T_{tP'u_1} \cap \ldots \cap  T_{tP'u_k}) >\frac{C}{4}.
\]

For $m=2$ the proof is similar.
By remark \ref{rmr:complex} there exists $c \in \BC$ 
such that for  all $t \in \BR$  we have
\[
\lim_{N \ra \infty} \Avr{N}{n} \p_{j=1}^k T_{ncu_j+itcu_j} E(f|Y_k)(y)=g(y)
\]
in $L^2(Y)$, and 
\[
\Avr{N}{n} \p_{j=1}^k T_{ncu_j+itcu_j} f(x)-
\Avr{N}{n} \p_{j=1}^k T_{ncu_j+itcu_j} \pi^*E(f|Y_k)(x)\ra 0
\]
in $L^2(X)$, where the convergence is uniform in $t$. 
As in the proof for $m>2$, there exists $N_0$, such that for all $t$
\[
\Avr{N_0}{n} \int f(x) \p_{j=1}^k T_{ncu_j+itcu_j} f(x) d\mu_X >\frac{C}{2}.
\]
This implies that for all  $t \in \BR$ there exists $n\le N_0$ 
with  
\begin{equation*}
\begin{split}
\mu( A\cap & T_{(n+it)cu_1}A\cap \ldots \cap T_{(n+it)cu_k}A) = \\
&\int f(x)\p_{j=1}^k T_{(n+it)cu_j} f(x) d\mu_X 
>\frac{C}{2}.
\end{split}
\end{equation*}
If $t$ is large enough, then 
\[
(n+it)cu_j \sim \frac{t}{|n+it|}(n+it)cu_j,
\]
and $|\frac{t}{|n+it|}(n+it)|=t$.
\end{dsc}

%%%%%%%%%%%%%%%%%%%%%
% Bibliography      %
%%%%%%%%%%%%%%%%%%%%%

\end{document}